\documentclass{amsart}
\usepackage{latexsym,xypic,amssymb,amscd,amsmath}
\usepackage{graphics}

\theoremstyle{plain}
\newtheorem{theorem}{Theorem}
\newtheorem{lemma}{Lemma}
\newtheorem{corollary}{Corollary}

\theoremstyle{definition}

\begin{document}
\title{Algebraic monoids with affine unit group are affine}
\author{Alvaro Rittatore} 

%
\maketitle

\begin{abstract}
In this short note we prove that any irreducible algebraic monoid whose unit
group is an affine algebraic group is affine. 
\end{abstract}

\section{Introduction}

Let $\Bbbk$ be an algebraically closed field of arbitrary
characteristic. An {\em  algebraic monoid}\/ is an 
algebraic variety  $M$ with an associative  product $M\times
M\rightarrow M$ which 
is a morphism of algebraic varieties, such that there exists a neutral
element $1$ for this 
product. In this case, it can be proved that the unit group $G(M)$
-- i.e.~the group of invertible elements -- is an 
algebraic group, open in $M$ (see \cite[Thm. 1]{kn:am} and Lemma \ref{lema:density}).  

It is easy to show that the action   $\bigl(G(M)\times
G(M)\bigr) \times M\rightarrow M$, $\bigl((a,b),m\bigr)\mapsto
amb^{-1}$ is regular, with open orbit $G(M)\cong \bigl( G(M)\times
G(M)\bigr)/\Delta\bigl(G(M)\bigr)$, where $\Delta\bigl(G(M)\bigr)$ is
the diagonal. In other words, $M$ is a {\em
  $G(M)$-embedding}.  Moreover, 
it is a {\em simple}\/ embedding, i.e.~there exists an unique closed
orbit, namely  the {\em center}\/ of $M$ -- the
minimum ideal of $M$ --  (see \cite[Thm. 1]{kn:am} and Lemma \ref{lema:density}).

It is well known that if $G$ is a quasi-affine algebraic group, then
$G$ is affine (see for example \cite[Thm. 7.5.3]{kn:libro}); 
in \cite[Thm.\ 4.4]{kn:renner-quasiaffine}, Renner proved the analog
for algebraic monoids, namely that
any quasi-affine algebraic monoid is affine. In particular, 
if $M$ is an affine algebraic monoid, then $G(M)$ is quasi-affine and
hence affine. Conversely, in \cite[Prop. 1]{kn:am} it is proved that
if $M$ is  an irreducible affine embedding of a (necesarily affine) algebraic
group $G$, then $M$ is an algebraic monoid of unit group $G(M)=G$. 
These observations lead naturally to the following conjecture
(presented as an open problem by Renner in \cite{kn:renner-quasiaffine}),
communicated to the author by E.B.~Vinberg:

{\em Let $M$ be an irreducible algebraic monoid whose unit group $G(M)$ is
  affine, then $M$ is also affine.}

A partial affirmative answer was given in \cite{kn:am}, where it is
proved that any irreducible reductive monoid -- i.e.\ with reductive
unit group -- is affine, and those monoids are classified in
combinatorial terms. 
In this note we give an affirmative answer to the above conjecture 
(see Theorem \ref{theo:elteorema}). Our methods are based on a
generalization to the context of algebraic
monoids of results by F.~Knop, H.~Kraft, D.~Luna and  T.~Vust about line
bundles  over affine algebraic groups (\cite{kn:KKLV}). In the last
section, we deal with the non-irreducible case, showing that any
algebraic monoid $M$ with unit group  affine and 
  dense in $M$ is an affine algebraic variety. 

The author would like to thank M.~Brion and W.\ Ferrer Santos for many
useful suggestions and 
remarks.

\section{Main results}

We begin this section by recalling a well known fact about algebraic
monoids (see \cite{kn:renner-quasiaffine}).

\begin{lemma}
\label{lema:monoidquasi-proj}
Let $M$ be a normal irreducible algebraic monoid with affine unit
group $G$. Then $M$ is quasi-projective. 
\end{lemma}

\proof
 Since $M$ is a simple $G$-embedding, of closed orbit its 
 center, we  can apply 
Sumihiro's theorem  (\cite{kn:sumihiro}) and obtain a finite dimensional 
$(G\times G)$-module $V$ and a $(G\times G)$-equivariant  open immersion
$M\hookrightarrow \mathbb P(V)$. 
\qed

Let $X$ be an algebraic variety, and $\pi:L\rightarrow X$ the line
bundle associated 
to an invertible coherent sheaf on $X$. Recall that the {\em zero
  section}\/  $\sigma_0 : X\rightarrow L$ is defined in the following
way: let $U\subset X$ be a trivializing  open subset of $L$, and identify
  $\pi^{-1}(U)\cong U\times \Bbbk$. Then 
  $\sigma_0|_{_U}:U\rightarrow U\times \Bbbk$, $\sigma_0(x)=(x,0)$
  (it is an easy exercise to prove that this is well defined, see
  for example \cite[p.\ 128]{kn:hartshorne}).    
We will denote as $L^*=L \setminus \sigma_0(X)$.

The following key result is due to Demazure and Fujita (see \cite[Lemma 1.1.13]{kn:brion-kumar},
\cite{kn:demazure-anneaux-gradues} or \cite{kn:fujita}).

\begin{lemma}
\label{lema:fujita}
Let $X$ be a projective normal algebraic variety, $\mathcal L$ an
ample invertible sheaf, and $\pi:L\rightarrow
X$ the line bundle associated to $\mathcal L^{-1}$, the dual of
$\mathcal L$. Then

\noindent (1) the algebra $R(X,L)=\bigoplus_{n\geq 0}H^0(X,L^{\otimes
n})$ is affine and integrally closed on its field of fractions; 

\noindent (2) the morphism $\varphi:L\rightarrow
\operatorname{Spec}R(X,L)$ induced by the restriction of sections is
proper with connected fibres;

\noindent (3) the restriction  $\varphi|_{_{L\setminus
\sigma_0(X)}}$ induces an isomorphism 
\[
L^*=L\setminus
\sigma_0(X) \cong\operatorname{Spec}R(X,L)\setminus \{O\},
\]
 where $O$
is the point associated to the maximal ideal $R^+(X,L)= {\displaystyle 
\bigoplus_{n\geq 1}}H^0(X,L^{\otimes n})$.\qed  
\end{lemma}

\noindent {\bf Remark} Property (1) holds even if only a power of
$\mathcal L$ is generated by  
its global sections (\cite{kn:fujita}). Demazure proves (3) in  \cite[\S
3]{kn:demazure-anneaux-gradues}, for Weil divisors with rational
coefficients such that a positive power is ample. Moreover, he shows
that  in this way one obtains all the finitely generated normal
graded algebras.

\medskip

The following result is a generalization of
\cite[Lemmas 4.2 and 4.3]{kn:KKLV}.

\begin{lemma}
\label{lema:knopgen}
Let $G$ be an affine algebraic group and $X$ a projective normal $G$-embedding, 
 $\pi:L\to X$ be  a $(G\times G)$-linearized 
line bundle, and  
consider  the commutative diagram:

\begin{center}
\mbox{
\xymatrix{
(i^*L)^* \ar@{->}[d]   \ar@{->}[r] & L\ar@{->}[d] \\ 
G \ar@{^(->}[r]_{i}   & X
}
}
\end{center}

If we denote as  $H=(i^*L)^*=(i^*L)\setminus 
\sigma_0(G)$, then $H\cong G\times \Bbbk^*$ and  $L^*$ is an $H$-embedding.
\end{lemma}

\proof
Since $L$ is $(G\times G)$-linearized, then $G$ acts on the fibre
$\pi^{-1}(1)\cong \Bbbk$ as the diagonal $\Delta(G)\subset G\times
G$, i.e.\ $g\cdot l=(g,g)\cdot l$ for all $l\in \pi^{-1}(1)$. Hence,
$G$ acts by multiplication with a character $\lambda:G\to \Bbbk^*$. 
Moreover, for all $g\in G$ and  $l\in \pi^{-1}(g)$ we have that 
\[
(a,g^{-1}ag)\cdot
l=\bigl((1,g^{-1})(a,a)(a,g)\bigr)\cdot l=
(1,g^{-1})\bigl(\lambda(a)(1,g)\cdot l\bigr) = \lambda(a) l.
\] 

Extend $\lambda:G\to \Bbbk^*$  to $G\times G$ by 
$\widetilde{\lambda}(g,g')
= \lambda(g)$, and change the linearization  by considering $(a,b)\star
l=\widetilde{\lambda}^{-1}(a,b)\cdot l=\lambda^{-1}(a)(a,b)\cdot
l$. Then $L$ is trivial over $G$, and $H=(i^*L)^*\cong G\times \Bbbk^*$ is an
algebraic group such that $H\times H$ acts on $L^*$ by
$\bigl((g,s),(g',t)\bigr)\cdot l= st^{-1}\bigl((g,g')\star
l\bigr)$; it is clear that $H$ is an open orbit for this action. 
\qed

\begin{theorem}
\label{theo:conemonoid}
Let $M$ be a normal irreducible algebraic monoid with unit group
 an affine algebraic group
$G$. Then there exists  a  $(G\times G)$-linearized line bundle
$\pi:N\to M$, such that  $N^*$ is an affine algebraic monoid with unit
group $H=\pi|_{_{N^*}}^{-1}(G)\cong G\times \Bbbk^*$. Moreover,
 $\pi:N^*\to M$ is a morphism of algebraic monoids, and is the quotient of
 $N^*$ by $\pi^{-1}(1)\cong \Bbbk^*$.
\end{theorem}

\proof
By Sumihiro's theorem, there exists an open $(G\times G)$-equi\-va\-riant
inmersion $\varphi :M\hookrightarrow X$, where $X$ is a projective
$G$-embedding. It is clear that we can suppose that $X$ is normal,
and that there exists a very  ample invertible $(G\times
G)$-linearizable sheaf $\mathcal L$ on $X$ (see \cite[Proposition
2.4]{kn:KKLV}). Let $\pi:L\rightarrow X$ be the line bundle associated
to the dual of $\mathcal L$,  and 
$N=\varphi^*(M)$ its  restriction to $M$.

By Lemma \ref{lema:knopgen}, $N$ is a
$H=\pi|_{_{N^*}}^{-1}(G)$-embedding. Let $\mu_1: H\times L^*\to L^*$,
$\mu_1(h,l)=(h,1)\cdot l$, $\mu_2:L^*\times H\to L^*$,
$\mu_2(l,h)=(1,h^{-1})\cdot l$. Since both  $\mu_1$ and $\mu_2$
coincide  with the product on $H$ when restricted to $H\times
H=(H\times L^*)\cap (L^*\times H)$, they induce a morphism
 $\mu: U=H\times L^*\cup
L^*\times H\rightarrow L^*$. Since  
$L^*=\operatorname{Spec}R(X,L)\setminus \{O\}$ is a quasi-affine
normal variety by Lemma \ref{lema:fujita}, 
and clearly $\operatorname{codim} (L^*\times L^*)\setminus U\geq 2$, 
we can extend the morphism $\mu$ to a morphism  $\mu:L^*\times
L^*\to \operatorname{Spec}R(X,L)$. 

It suffices to prove that $\mu(N^*\times N^*)\subset N^*$. Indeed,
if this is the case then $\mu$ is an associative product in $N^*$,
since it is associative on the open subset $H\times H\times H\subset
N^*\times N^*\times N^*$. Then $N^*$ is a quasi-affine algebraic
monoid,  and  it follows from the result of Renner cited at the
introduction (see \cite[Thm. 4.4]{kn:renner-quasiaffine}) that $N^*$
 is an affine algebraic monoid.  

In order to prove that $\mu(N^*\times N^*)\subset N^*$, consider $u,v\in
N^*$. There exists an affine open subset $V\subset M$ such that
 $\pi(u)\pi(v)\in V$ and $\pi^{-1}(V)\cong V\times \Bbbk^*$. 
Let $m:M\times  M\to M$ be the product and consider
$W=m^{-1}(V)\subset M\times M$; let $W'=\pi^{-1}(W)\cap 
\bigl( (H\times N^*)\cup (N^*\times H)\bigr)$. Then $W'$ is
an open subset of $\pi^{-1}(W)$,   with complement of codimesion 
greater than $2$, such that $\mu(W')\subset \pi^{-1}(V)$. Since $N^*$
and hence $\pi^{-1}(W)$ are normal, it follows 
that $\mu|_{_{W'}}:W'\to  \pi^{-1}(V)$ extends to a morphism
$\widetilde{\mu}:\pi^{-1}(W)\to  \pi^{-1}(V)$. Since both $\mu$ and
$\widetilde{\mu}$ are continous functions and
$\mu|_{_{W'}}=\widetilde{\mu}|_{_{W'}}$,  it 
follows that $\mu|_{_{\pi^{-1}(W)}}=\widetilde{\mu}$; in particular,
$\mu(u,v)\in \pi^{-1}(V)\subset N^*$.

By construction, the map $\pi:N^*\to M$ is a morphism of algebraic
monoids, with central kernel  $\pi^{-1}(1)\cong \Bbbk^*$. Hence, 
$M$ is a quotient of $N^*$ by
$\Bbbk^*$.
\qed

\begin{theorem}
\label{theo:elteorema}
Let $G$ be an affine algebraic group and $M$ an   algebraic
monoid with unit group $G$, affine algebraic group. Then $M$ is affine.
\end{theorem}

\proof
We can assume without loss of generality that $M$ is normal (see
for example \cite[Lemma 1]{kn:am}). Applying Theorem
 \ref{theo:conemonoid} we deduce that $M$
 is the quotient of an affine algebraic variety by an algebraic
 torus, and hence it is affine.
\qed

\section{The non-irreducible case}

Let $G$ be a non-connected affine algebraic group, and assume that $M$ is 
an algebraic monoid with unit group $G$. In order to obtain a better
control of the geometry of $M$ it is natural to impose the 
{\em density condition}\/  $\overline{G}=M$, as the following examples show:

\medskip

\noindent {\bf Examples}
(1) Let $S$ be an arbitrary algebraic variety and $s_0\in S$. 
Then $m:S\times S \to S$, $m(s,t)=s_0$  is an associative product. If
$M$ is an arbitrary  algebraic monoid, then the products on $M$ and
$S$ extend to a product $\mu$ on $M\cup S$   (disjoint union) by
$\mu(a,s)=s$ for all $a\in M$ and $s\in S$. Then $M\cup S$ is  
an algebraic monoid, of unit group $G(M)$ and zero $s_0$. 

\noindent (2) Assume now that $M$ has a zero $0$,  and consider 
the equivalence relationship on $M\cup S$ induced by $0\sim s_0$. 
Then $\mu$ induces a product $\widetilde{\mu}$ on $N=(M\cup S)/{\sim}$, 
in such a way that $N$ is an algebraic monoid with unit group $G(M)$. 
Observe that $N$  can be realized as the closed subvariety 
 $N\cong\bigl(M\times \{s_0\}\bigr)\cup
\bigl(\{0\}\times S\bigr) \subset M\times S$.

The following lemma is an easy generalization of \cite[Thm.\
1]{kn:am}, where  the case of 
irreducible algebraic monoids is treated, hence we omit the proof.

\begin{lemma}
\label{lema:density}
Let $M$ be an algebraic monoid of unit group $G$, dense in $M$. Then
$M$ is a simple  
$G$-embedding, with unique closed orbit the center of $M$.
\qed
\end{lemma}

The following theorem generalizes \cite[Prop.\ 2]{kn:am},
where  it is proved that if an algebraic monoid
verifies the density condition $\overline{G(M)}=M$, then any two
irreducible components are 
isomorphic.

\begin{theorem}
Let $M$ be an algebraic monoid with affine dense unit group $G$ and center
$Y$. Let  $G=\cup_{i=1}^nG_i$, where $1\in G_1$, be the decomposition in
irreducible  components of $G$; then  $M=\cup_{i=1}^n 
\overline{G_i}$. If we set $M_i=\overline{G_i}$, $i=1,\dots ,n$, then
$M_1$ is an affine algebraic monoid of unit group $G_1$,
and $M_i\cong M_1$ as an algebraic variety for all $i=1,\dots, n$. In particular, $M$ is an
affine algebraic variety. 

Moreover, $M_i\cap Y \subset G_i\cdot Y_1=Y_1\cdot G_i$, where $Y_1$ denotes the center of $M_1$.
\end{theorem}

\proof
Consider  $\mu: M\times M\to M$, the product on $M$; then for all
$j=1,\dots, n$, $M_1\cdot M_j= \mu(M_1\times M_j)$ is irreducible and
contains $1\cdot M_j=M_j$. Hence $M_1\cdot M_j=M_j$; in particular,
$M_1$ is an algebraic monoid with unit group $G_1$ and it follows from
Theorem \ref{theo:elteorema} that $M_1$ is an  affine algebraic
variety. Moreover, if $g_i\in G_i$, then $\ell_{g_i}:M_1\to M_i$,
$\ell_{g_i}(m)=g_im$   is an isomorphism with inverse $\ell_{g_i^{-1}}$,
and thus $M$ is an affine algebraic variety.

In order to prove the last assertion, observe that $G_i\cdot
 Y_1\subset M_i\cdot Y_1\subset M_i\cap Y$, and that since $G_1$ is
 normal in $G$, it follows that $G_i\cdot Y_1=Y_1\cdot G_i$.
\qed 

\begin{corollary}
Let $M$ be an algebraic monoid with zero $0\in M$ satisfying the density
condition. Then every irreducible component of $M$ contains $0$.
\qed
\end{corollary}

\bibliographystyle{amsplain}

\bibliography{oam-biblio}

\bigskip

\noindent {\sc Alvaro Rittatore}\\
Facultad de Ciencias\\ 
Universidad de la Rep\'ublica\\
Igu\'a 4225\\
11400 Montevideo\\
Uruguay\\
e-mail: {\tt alvaro@cmat.edu.uy}

\end{document}